

\ifx\shlhetal\undefinedcontrolsequence\let\shlhetal\relax\fi

\input amstex
\expandafter\ifx\csname mathdefs.tex\endcsname\relax
  \expandafter\gdef\csname mathdefs.tex\endcsname{}
\else \message{Hey!  Apparently you were trying to
  \string\input{mathdefs.tex} twice.   This does not make sense.} 
\errmessage{Please edit your file (probably \jobname.tex) and remove
any duplicate ``\string\input'' lines}\endinput\fi




\catcode`\X=12\catcode`\@=11

\def\n@wcount{\alloc@0\count\countdef\insc@unt}
\def\n@wwrite{\alloc@7\write\chardef\sixt@@n}
\def\n@wread{\alloc@6\read\chardef\sixt@@n}
\def\r@s@t{\relax}\def\v@idline{\par}\def\@mputate#1/{#1}
\def\l@c@l#1X{\firstpart.#1}\def\gl@b@l#1X{#1}\def\t@d@l#1X{{}}

\def\crossrefs#1{\ifx\all#1\let\tr@ce=\all\else\def\tr@ce{#1,}\fi
   \n@wwrite\cit@tionsout\openout\cit@tionsout=\jobname.cit 
   \write\cit@tionsout{\tr@ce}\expandafter\setfl@gs\tr@ce,}
\def\setfl@gs#1,{\def\@{#1}\ifx\@\empty\let\next=\relax
   \else\let\next=\setfl@gs\expandafter\xdef
   \csname#1tr@cetrue\endcsname{}\fi\next}
\def\m@ketag#1#2{\expandafter\n@wcount\csname#2tagno\endcsname
     \csname#2tagno\endcsname=0\let\tail=\all\xdef\all{\tail#2,}
   \ifx#1\l@c@l\let\tail=\r@s@t\xdef\r@s@t{\csname#2tagno\endcsname=0\tail}\fi
   \expandafter\gdef\csname#2cite\endcsname##1{\expandafter
     \ifx\csname#2tag##1\endcsname\relax?\else\csname#2tag##1\endcsname\fi
     \expandafter\ifx\csname#2tr@cetrue\endcsname\relax\else
     \write\cit@tionsout{#2tag ##1 cited on page \folio.}\fi}
   \expandafter\gdef\csname#2page\endcsname##1{\expandafter
     \ifx\csname#2page##1\endcsname\relax?\else\csname#2page##1\endcsname\fi
     \expandafter\ifx\csname#2tr@cetrue\endcsname\relax\else
     \write\cit@tionsout{#2tag ##1 cited on page \folio.}\fi}
   \expandafter\gdef\csname#2tag\endcsname##1{\expandafter
      \ifx\csname#2check##1\endcsname\relax
      \expandafter\xdef\csname#2check##1\endcsname{}%
      \else\immediate\write16{Warning: #2tag ##1 used more than once.}\fi
      \multit@g{#1}{#2}##1/X%
      \write\t@gsout{#2tag ##1 assigned number \csname#2tag##1\endcsname\space
      on page \number\count0.}%
   \csname#2tag##1\endcsname}}

\def\multit@g#1#2#3/#4X{\def\t@mp{#4}\ifx\t@mp\empty%
      \global\advance\csname#2tagno\endcsname by 1 
      \expandafter\xdef\csname#2tag#3\endcsname
      {#1\number\csname#2tagno\endcsnameX}%
   \else\expandafter\ifx\csname#2last#3\endcsname\relax
      \expandafter\n@wcount\csname#2last#3\endcsname
      \global\advance\csname#2tagno\endcsname by 1 
      \expandafter\xdef\csname#2tag#3\endcsname
      {#1\number\csname#2tagno\endcsnameX}
      \write\t@gsout{#2tag #3 assigned number \csname#2tag#3\endcsname\space
      on page \number\count0.}\fi
   \global\advance\csname#2last#3\endcsname by 1
   \def\t@mp{\expandafter\xdef\csname#2tag#3/}%
   \expandafter\t@mp\@mputate#4\endcsname
   {\csname#2tag#3\endcsname\lastpart{\csname#2last#3\endcsname}}\fi}
\def\t@gs#1{\def\all{}\m@ketag#1e\m@ketag#1s\m@ketag\t@d@l p
\let\realscite\scite
\let\realstag\stag
   \m@ketag\gl@b@l r \n@wread\t@gsin
   \openin\t@gsin=\jobname.tgs \re@der \closein\t@gsin
   \n@wwrite\t@gsout\openout\t@gsout=\jobname.tgs }
\outer\def\localtags{\t@gs\l@c@l}
\outer\def\globaltags{\t@gs\gl@b@l}
\outer\def\newlocaltag#1{\m@ketag\l@c@l{#1}}
\outer\def\newglobaltag#1{\m@ketag\gl@b@l{#1}}

\newif\ifpr@ 
\def\m@kecs #1tag #2 assigned number #3 on page #4.%
   {\expandafter\gdef\csname#1tag#2\endcsname{#3}
   \expandafter\gdef\csname#1page#2\endcsname{#4}
   \ifpr@\expandafter\xdef\csname#1check#2\endcsname{}\fi}
\def\re@der{\ifeof\t@gsin\let\next=\relax\else
   \read\t@gsin to\t@gline\ifx\t@gline\v@idline\else
   \expandafter\m@kecs \t@gline\fi\let \next=\re@der\fi\next}
\def\pretags#1{\pr@true\pret@gs#1,,}
\def\pret@gs#1,{\def\@{#1}\ifx\@\empty\let\n@xtfile=\relax
   \else\let\n@xtfile=\pret@gs \openin\t@gsin=#1.tgs \message{#1} \re@der 
   \closein\t@gsin\fi \n@xtfile}

\newcount\sectno\sectno=0\newcount\subsectno\subsectno=0
\newif\ifultr@local \def\ultralocal{\ultr@localtrue}
\def\firstpart{\number\sectno}
\def\lastpart#1{\ifcase#1 \or a\or b\or c\or d\or e\or f\or g\or h\or 
   i\or k\or l\or m\or n\or o\or p\or q\or r\or s\or t\or u\or v\or w\or 
   x\or y\or z \fi}

\def\resetall{\global\advance\sectno by 1\subsectno=0
   \gdef\firstpart{\number\sectno}\r@s@t}
\def\resetsub{\global\advance\subsectno by 1
   \gdef\firstpart{\number\sectno.\number\subsectno}\r@s@t}
\def\newsection#1\par{\resetall\vskip0pt plus.3\vsize\penalty-250
   \vskip0pt plus-.3\vsize\bigskip\bigskip
   \message{#1}\leftline{\bf#1}\nobreak\bigskip}
\def\subsection#1\par{\ifultr@local\resetsub\fi
   \vskip0pt plus.2\vsize\penalty-250\vskip0pt plus-.2\vsize
   \bigskip\smallskip\message{#1}\leftline{\bf#1}\nobreak\medskip}


\newdimen\marginshift

\newdimen\margindelta
\newdimen\marginmax
\newdimen\marginmin

\def\margininit{       
\marginmax=3 true cm                  
				      
\margindelta=0.1 true cm              
\marginmin=0.1true cm                 
\marginshift=\marginmin
}    

\def\t@gsjj#1,{\def\@{#1}\ifx\@\empty\let\next=\relax\else\let\next=\t@gsjj
   \def\@@{p}\ifx\@\@@\else
   \expandafter\gdef\csname#1cite\endcsname##1{\citejj{##1}}
   \expandafter\gdef\csname#1page\endcsname##1{?}
   \expandafter\gdef\csname#1tag\endcsname##1{\tagjj{##1}}\fi\fi\next}
\newif\ifshowstuffinmargin
\showstuffinmarginfalse
\def\jjtags{\ifx\shlhetal\relax 
  \else
\ifx\shlhetal\undefinedcontrolseq
\else
\showstuffinmargintrue
\ifx\all\relax\else\expandafter\t@gsjj\all,\fi\fi \fi
}

\def\tagjj#1{\realstag{#1}\oldmginpar{\zeigen{#1}}}
\def\citejj#1{\rechnen{#1}\mginpar{\zeigen{#1}}}     

\def\rechnen#1{\expandafter\ifx\csname stag#1\endcsname\relax ??\else
                           \csname stag#1\endcsname\fi}

\newdimen\theight

\def\marginfont{\sevenrm}

\def\trymarginbox#1{\setbox0=\hbox{\marginfont\hskip\marginshift #1}%
		\global\marginshift\wd0 
		\global\advance\marginshift\margindelta}

\def \oldmginpar#1{%
\ifvmode\setbox0\hbox to \hsize{\hfill\rlap{\marginfont\quad#1}}%
\ht0 0cm
\dp0 0cm
\box0\vskip-\baselineskip
\else 
             \vadjust{\trymarginbox{#1}%
		\ifdim\marginshift>\marginmax \global\marginshift\marginmin
			\trymarginbox{#1}%
                \fi
             \theight=\ht0
             \advance\theight by \dp0    \advance\theight by \lineskip
             \kern -\theight \vbox to \theight{\rightline{\rlap{\box0}}%
\vss}}\fi}

\newdimen\upordown
\global\upordown=8pt
\font\tinyfont=cmtt8 
\def\mginpar#1{\smash{\hbox to 0cm{\kern-10pt\raise7pt\hbox{\tinyfont #1}\hss}}}
\def\mginpar#1{{\hbox to 0cm{\kern-10pt\raise\upordown\hbox{\tinyfont #1}\hss}}\global\upordown-\upordown}


\def\t@gsoff#1,{\def\@{#1}\ifx\@\empty\let\next=\relax\else\let\next=\t@gsoff
   \def\@@{p}\ifx\@\@@\else
   \expandafter\gdef\csname#1cite\endcsname##1{\zeigen{##1}}
   \expandafter\gdef\csname#1page\endcsname##1{?}
   \expandafter\gdef\csname#1tag\endcsname##1{\zeigen{##1}}\fi\fi\next}
\def\verbatimtags{\showstuffinmarginfalse
\ifx\all\relax\else\expandafter\t@gsoff\all,\fi}
\def\zeigen#1{\hbox{$\scriptstyle\langle$}#1\hbox{$\scriptstyle\rangle$}}


\def\margintag#1{\ifshowstuffinmargin\oldmginpar{\zeigen{#1}}\fi}

\def\(#1){\edef\dot@g{\ifmmode\ifinner(\hbox{\noexpand\etag{#1}})
   \else\noexpand\eqno(\hbox{\noexpand\etag{#1}})\fi
   \else(\noexpand\ecite{#1})\fi}\dot@g}

\newif\ifbr@ck
\def\eat#1{}
\def\[#1]{\br@cktrue[\br@cket#1'X]}
\def\br@cket#1'#2X{\def\temp{#2}\ifx\temp\empty\let\next\eat
   \else\let\next\br@cket\fi
   \ifbr@ck\br@ckfalse\br@ck@t#1,X\else\br@cktrue#1\fi\next#2X}
\def\br@ck@t#1,#2X{\def\temp{#2}\ifx\temp\empty\let\neext\eat
   \else\let\neext\br@ck@t\def\temp{,}\fi
   \def\teemp{#1}\ifx\teemp\empty\else\rcite{#1}\fi\temp\neext#2X}
\def\resetbr@cket{\gdef\[##1]{[\rtag{##1}]}}
\def\references{\resetbr@cket\newsection References\par}

\newtoks\symb@ls\newtoks\s@mb@ls\newtoks\p@gelist\n@wcount\ftn@mber
    \ftn@mber=1\newif\ifftn@mbers\ftn@mbersfalse\newif\ifbyp@ge\byp@gefalse
\def\defm@rk{\ifftn@mbers\n@mberm@rk\else\symb@lm@rk\fi}
\def\n@mberm@rk{\xdef\m@rk{{\the\ftn@mber}}%
    \global\advance\ftn@mber by 1 }
\def\rot@te#1{\let\temp=#1\global#1=\expandafter\r@t@te\the\temp,X}
\def\r@t@te#1,#2X{{#2#1}\xdef\m@rk{{#1}}}
\def\b@@st#1{{$^{#1}$}}\def\str@p#1{#1}
\def\symb@lm@rk{\ifbyp@ge\rot@te\p@gelist\ifnum\expandafter\str@p\m@rk=1 
    \s@mb@ls=\symb@ls\fi\write\f@nsout{\number\count0}\fi \rot@te\s@mb@ls}
\def\byp@ge{\byp@getrue\n@wwrite\f@nsin\openin\f@nsin=\jobname.fns 
    \n@wcount\currentp@ge\currentp@ge=0\p@gelist={0}
    \re@dfns\closein\f@nsin\rot@te\p@gelist
    \n@wread\f@nsout\openout\f@nsout=\jobname.fns }
\def\m@kelist#1X#2{{#1,#2}}
\def\re@dfns{\ifeof\f@nsin\let\next=\relax\else\read\f@nsin to \f@nline
    \ifx\f@nline\v@idline\else\let\t@mplist=\p@gelist
    \ifnum\currentp@ge=\f@nline
    \global\p@gelist=\expandafter\m@kelist\the\t@mplistX0
    \else\currentp@ge=\f@nline
    \global\p@gelist=\expandafter\m@kelist\the\t@mplistX1\fi\fi
    \let\next=\re@dfns\fi\next}
\def\symbols#1{\symb@ls={#1}\s@mb@ls=\symb@ls} 
\def\bigsymbol{\textstyle}
\symbols{\bigsymbol\ast,\dagger,\ddagger,\sharp,\flat,\natural,\star}
\def\ftnumbers{\ftn@mberstrue} \def\ftsymbols{\ftn@mbersfalse}
\def\paginal{\byp@ge} \def\resetftnumbers{\ftn@mber=1}
\def\ftnote#1{\defm@rk\expandafter\expandafter\expandafter\footnote
    \expandafter\b@@st\m@rk{#1}}

\long\def\jump#1\endjump{}
\def\ssum{\mathop{\lower .1em\hbox{$\textstyle\Sigma$}}\nolimits}

\def\qed{\nobreak\kern 1em \vrule height .5em width .5em depth 0em}
\def\newneq{\hbox{\rlap{\hbox to 1\wd9{\hss$=$\hss}}\raise .1em 
   \hbox to 1\wd9{\hss$\scriptscriptstyle/$\hss}}}
\def\subsetne{\setbox9 = \hbox{$\subset$}\mathrel{\hbox{\rlap
   {\lower .4em \newneq}\raise .13em \hbox{$\subset$}}}}
\def\supsetne{\setbox9 = \hbox{$\subset$}\mathrel{\hbox{\rlap
   {\lower .4em \newneq}\raise .13em \hbox{$\supset$}}}}

\def\vbar{\mathchoice{\vrule height6.3ptdepth-.5ptwidth.8pt\kern-.8pt}
   {\vrule height6.3ptdepth-.5ptwidth.8pt\kern-.8pt}
   {\vrule height4.1ptdepth-.35ptwidth.6pt\kern-.6pt}
   {\vrule height3.1ptdepth-.25ptwidth.5pt\kern-.5pt}}
\def\f@dge{\mathchoice{}{}{\mkern.5mu}{\mkern.8mu}}
\def\b@c#1#2{{\rm \mkern#2mu\vbar\mkern-#2mu#1}}
\def\b@b#1{{\rm I\mkern-3.5mu #1}}
\def\b@a#1#2{{\rm #1\mkern-#2mu\f@dge #1}}
\def\bb#1{{\count4=`#1 \advance\count4by-64 \ifcase\count4\or\b@a A{11.5}\or
   \b@b B\or\b@c C{5}\or\b@b D\or\b@b E\or\b@b F \or\b@c G{5}\or\b@b H\or
   \b@b I\or\b@c J{3}\or\b@b K\or\b@b L \or\b@b M\or\b@b N\or\b@c O{5} \or
   \b@b P\or\b@c Q{5}\or\b@b R\or\b@a S{8}\or\b@a T{10.5}\or\b@c U{5}\or
   \b@a V{12}\or\b@a W{16.5}\or\b@a X{11}\or\b@a Y{11.7}\or\b@a Z{7.5}\fi}}

\catcode`\X=11 \catcode`\@=12




\let\thischap\jobname

\def\partof#1{\csname returnthe#1part\endcsname}
\def\chapof#1{\csname returnthe#1chap\endcsname}

\def\setchapter#1,#2,#3;{%
  \expandafter\def\csname returnthe#1part\endcsname{#2}%
  \expandafter\def\csname returnthe#1chap\endcsname{#3}%
}

\setchapter 300a,A,II.A;
\setchapter 300b,A,II.B;
\setchapter 300c,A,II.C;
\setchapter 300d,A,II.D;
\setchapter 300e,A,II.E;
\setchapter 300f,A,II.F;
\setchapter 300g,A,II.G;
\setchapter  E53,B,N;
\setchapter  88r,B,I;
\setchapter  600,B,III;
\setchapter  705,B,IV;
\setchapter  734,B,V;

\def\cprefix#1{
\edef\theotherpart{\partof{#1}}\edef\theotherchap{\chapof{#1}}%
\ifx\theotherpart\thispart
   \ifx\theotherchap\thischap 
    \else 
     \theotherchap%
    \fi
   \else 
     \theotherchap\fi}

\def\sectioncite[#1]#2{%
     \cprefix{#2}#1}

\edef\thispart{\partof{\thischap}}
\edef\thischap{\chapof{\thischap}}

\def\lastpage of '#1' is #2.{\expandafter\def\csname lastpage#1\endcsname{#2}}



\expandafter\ifx\csname citeadd.tex\endcsname\relax
\expandafter\gdef\csname citeadd.tex\endcsname{}
\else \message{Hey!  Apparently you were trying to
\string\input{citeadd.tex} twice.   This does not make sense.} 
\errmessage{Please edit your file (probably \jobname.tex) and remove
any duplicate ``\string\input'' lines}\endinput\fi

\sectno=-1   
\localtags
\NoBlackBoxes
\define\mr{\medskip\roster}
\define\sn{\smallskip\noindent}
\define\mn{\medskip\noindent}
\define\bn{\bigskip\noindent}
\define\ub{\underbar}
\define\wilog{\text{without loss of generality}}
\define\ermn{\endroster\medskip\noindent}
\define\dbca{\dsize\bigcap}
\define\dbcu{\dsize\bigcup}
\define \nl{\newline}
\magnification=\magstep 1
\documentstyle{amsppt} 

{    
\catcode`@11

\ifx\alicetwothousandloaded@\relax
  \endinput\else\global\let\alicetwothousandloaded@\relax\fi

\gdef\subjclass{\let\savedef@\subjclass
 \def\subjclass##1\endsubjclass{\let\subjclass\savedef@
   \toks@{\def\usualspace{{\rm\enspace}}\eightpoint}%
   \toks@@{##1\unskip.}%
   \edef\thesubjclass@{\the\toks@
     \frills@{{\noexpand\rm2000 {\noexpand\it Mathematics Subject
       Classification}.\noexpand\enspace}}%
     \the\toks@@}}%
  \nofrillscheck\subjclass}
} 

\pageheight{8.5truein}
\topmatter
\title{On $T_3$-topological space omitting many cardinals} \endtitle
\author {Saharon Shelah \thanks {\null\newline 
This research was supported by the Israel Science Foundation and
I would like to thank Alice Leonhardt for the beautiful typing. 
Publ. No. 606} \endthanks} \endauthor 

\affil{Institute of Mathematics\\
 The Hebrew University\\
 Jerusalem, Israel
 \medskip
 Rutgers University\\
 Mathematics Department\\
 New Brunswick, NJ  USA} \endaffil

\abstract  We prove that for every (infinite cardinal) $\lambda$ there is
a $T_3$-space $X$ with clopen basis, $2^\mu$ points where $\mu = 2^\lambda$,
such that every closed subspace of cardinality $< |X|$ has cardinality 
$< \lambda$.
\endabstract
\endtopmatter
\document  

\expandafter\ifx\csname alice2jlem.tex\endcsname\relax
  \expandafter\xdef\csname alice2jlem.tex\endcsname{\the\catcode`@}
\else \message{Hey!  Apparently you were trying to
\string\input{alice2jlem.tex}  twice.   This does not make sense.}
\errmessage{Please edit your file (probably \jobname.tex) and remove
any duplicate ``\string\input'' lines}\endinput\fi

\expandafter\ifx\csname bib4plain.tex\endcsname\relax
  \expandafter\gdef\csname bib4plain.tex\endcsname{}
\else \message{Hey!  Apparently you were trying to \string\input
  bib4plain.tex twice.   This does not make sense.}
\errmessage{Please edit your file (probably \jobname.tex) and remove
any duplicate ``\string\input'' lines}\endinput\fi

\def\renewcommand{\newcommand}	       
\edef\cite{\the\catcode`@}%
\catcode`@ = 11
\let\@oldatcatcode = \cite
\chardef\@letter = 11
\chardef\@other = 12
%
%
%
%
\def\@innerdef#1#2{\edef#1{\expandafter\noexpand\csname #2\endcsname}}%
%
%
\@innerdef\@innernewcount{newcount}%
\@innerdef\@innernewdimen{newdimen}%
\@innerdef\@innernewif{newif}%
\@innerdef\@innernewwrite{newwrite}%
%
%
%
\def\@gobble#1{}%
%
%
%
\ifx\inputlineno\@undefined
   \let\@linenumber = \empty 
\else
   \def\@linenumber{\the\inputlineno:\space}%
\fi
%
%
%
\def\@futurenonspacelet#1{\def\cs{#1}%
   \afterassignment\@stepone\let\@nexttoken=
}%
\begingroup 
\def\\{\global\let\@stoken= }%
\\ 
\endgroup
\def\@stepone{\expandafter\futurelet\cs\@steptwo}%
\def\@steptwo{\expandafter\ifx\cs\@stoken\let\@@next=\@stepthree
   \else\let\@@next=\@nexttoken\fi \@@next}%
\def\@stepthree{\afterassignment\@stepone\let\@@next= }%
%
%
%
\def\@getoptionalarg#1{%
   \let\@optionaltemp = #1%
   \let\@optionalnext = \relax
   \@futurenonspacelet\@optionalnext\@bracketcheck
}%
%
%
\def\@bracketcheck{%
   \ifx [\@optionalnext
      \expandafter\@@getoptionalarg
   \else
      \let\@optionalarg = \empty
      \expandafter\@optionaltemp
   \fi
}%
\def\@@getoptionalarg[#1]{%
   \def\@optionalarg{#1}%
   \@optionaltemp
}%
%
%
%
\def\@nnil{\@nil}%
\def\@fornoop#1\@@#2#3{}%
\def\@for#1:=#2\do#3{%
   \edef\@fortmp{#2}%
   \ifx\@fortmp\empty \else
      \expandafter\@forloop#2,\@nil,\@nil\@@#1{#3}%
   \fi
}%
\def\@forloop#1,#2,#3\@@#4#5{\def#4{#1}\ifx #4\@nnil \else
       #5\def#4{#2}\ifx #4\@nnil \else#5\@iforloop #3\@@#4{#5}\fi\fi
}%
\def\@iforloop#1,#2\@@#3#4{\def#3{#1}\ifx #3\@nnil
       \let\@nextwhile=\@fornoop \else
      #4\relax\let\@nextwhile=\@iforloop\fi\@nextwhile#2\@@#3{#4}%
}%
%
%
%
\@innernewif\if@fileexists
\def\@testfileexistence{\@getoptionalarg\@finishtestfileexistence}%
\def\@finishtestfileexistence#1{%
   \begingroup
      \def\extension{#1}%
      \immediate\openin0 =
         \ifx\@optionalarg\empty\jobname\else\@optionalarg\fi
         \ifx\extension\empty \else .#1\fi
         \space
      \ifeof 0
         \global\@fileexistsfalse
      \else
         \global\@fileexiststrue
      \fi
      \immediate\closein0
   \endgroup
}%
%
%
%
%
\def\bibliographystyle#1{%
   \@readauxfile
   \@writeaux{\string\bibstyle{#1}}%
}%
\let\bibstyle = \@gobble
%
%
\let\bblfilebasename = \jobname
\def\bibliography#1{%
   \@readauxfile
   \@writeaux{\string\bibdata{#1}}%
   \@testfileexistence[\bblfilebasename]{bbl}%
   \if@fileexists
      \nobreak
      \@readbblfile
   \fi
}%
\let\bibdata = \@gobble
%
%
\def\nocite#1{%
   \@readauxfile
   \@writeaux{\string\citation{#1}}%
}%
\@innernewif\if@notfirstcitation
%
%
\def\cite{\@getoptionalarg\@cite}%
%
%
\def\@cite#1{%
   \let\@citenotetext = \@optionalarg
   \printcitestart
   \nocite{#1}%
   \@notfirstcitationfalse
   \@for \@citation :=#1\do
   {%
      \expandafter\@onecitation\@citation\@@
   }%
   \ifx\empty\@citenotetext\else
      \printcitenote{\@citenotetext}%
   \fi
   \printcitefinish
}%
\newif\ifweareinprivate
\weareinprivatetrue
\ifx\shlhetal\undefinedcontrolseq\weareinprivatefalse\fi
\ifx\shlhetal\relax\weareinprivatefalse\fi
\def\@onecitation#1\@@{%
   \if@notfirstcitation
      \printbetweencitations
   \fi
   \expandafter \ifx \csname\@citelabel{#1}\endcsname \relax
      \if@citewarning
         \message{\@linenumber Undefined citation `#1'.}%
      \fi
     \ifweareinprivate
      \expandafter\gdef\csname\@citelabel{#1}\endcsname{%
\strut 
\vadjust{\vskip-\dp\strutbox
\vbox to 0pt{\vss\parindent0cm \leftskip=\hsize 
\advance\leftskip3mm
\advance\hsize 4cm\strut\openup-4pt 
\rightskip 0cm plus 1cm minus 0.5cm ?  #1 ?\strut}}
         {\tt
            \escapechar = -1
            \nobreak\hskip0pt\pfeilsw
            \expandafter\string\csname#1\endcsname
             \pfeilso
            \nobreak\hskip0pt
         }%
      }%
     \else  
      \expandafter\gdef\csname\@citelabel{#1}\endcsname{%
            {\tt\expandafter\string\csname#1\endcsname}
      }%
     \fi  
   \fi
   \csname\@citelabel{#1}\endcsname
   \@notfirstcitationtrue
}%
%
%
\def\@citelabel#1{b@#1}%
%
%
\def\@citedef#1#2{\expandafter\gdef\csname\@citelabel{#1}\endcsname{#2}}%
%
%
%
\def\@readbblfile{%
   \ifx\@itemnum\@undefined
      \@innernewcount\@itemnum
   \fi
   \begingroup
      \def\begin##1##2{%
         \setbox0 = \hbox{\biblabelcontents{##2}}%
         \biblabelwidth = \wd0
      }%
      \def\end##1{}
      %
      %
      \@itemnum = 0
      \def\bibitem{\@getoptionalarg\@bibitem}%
      \def\@bibitem{%
         \ifx\@optionalarg\empty
            \expandafter\@numberedbibitem
         \else
            \expandafter\@alphabibitem
         \fi
      }%
      \def\@alphabibitem##1{%
         \expandafter \xdef\csname\@citelabel{##1}\endcsname {\@optionalarg}%
         \ifx\biblabelprecontents\@undefined
            \let\biblabelprecontents = \relax
         \fi
         \ifx\biblabelpostcontents\@undefined
            \let\biblabelpostcontents = \hss
         \fi
         \@finishbibitem{##1}%
      }%
      \def\@numberedbibitem##1{%
         \advance\@itemnum by 1
         \expandafter \xdef\csname\@citelabel{##1}\endcsname{\number\@itemnum}%
         \ifx\biblabelprecontents\@undefined
            \let\biblabelprecontents = \hss
         \fi
         \ifx\biblabelpostcontents\@undefined
            \let\biblabelpostcontents = \relax
         \fi
         \@finishbibitem{##1}%
      }%
      \def\@finishbibitem##1{%
         \biblabelprint{\csname\@citelabel{##1}\endcsname}%
         \@writeaux{\string\@citedef{##1}{\csname\@citelabel{##1}\endcsname}}%
         \ignorespaces
      }%
      %
      %
      \let\em = \bblem
      \let\newblock = \bblnewblock
      \let\sc = \bblsc
      \frenchspacing
      \clubpenalty = 4000 \widowpenalty = 4000
      \tolerance = 10000 \hfuzz = .5pt
      \everypar = {\hangindent = \biblabelwidth
                      \advance\hangindent by \biblabelextraspace}%
      \bblrm
      \parskip = 1.5ex plus .5ex minus .5ex
      \biblabelextraspace = .5em
      \bblhook
      \input \bblfilebasename.bbl
   \endgroup
}%
%
%
\@innernewdimen\biblabelwidth
\@innernewdimen\biblabelextraspace
%
%
%
\def\biblabelprint#1{%
   \noindent
   \hbox to \biblabelwidth{%
      \biblabelprecontents
      \biblabelcontents{#1}%
      \biblabelpostcontents
   }%
   \kern\biblabelextraspace
}%
%
%
%
\def\biblabelcontents#1{{\bblrm [#1]}}%
%
%
\def\bblrm{\rm}%
%
%
\def\bblem{\it}%
%
%
\def\bblsc{\ifx\@scfont\@undefined
              \font\@scfont = cmcsc10
           \fi
           \@scfont
}%
%
%
\def\bblnewblock{\hskip .11em plus .33em minus .07em }%
%
%
\let\bblhook = \empty
%
%
%
\def\printcitestart{[}
\def\printcitefinish{]}
\def\printbetweencitations{, }
\def\printcitenote#1{, #1}
%
%
%
\let\citation = \@gobble
%
%
%
\@innernewcount\@numparams
%
%
\def\newcommand#1{%
   \def\@commandname{#1}%
   \@getoptionalarg\@continuenewcommand
}%
%
%
\def\@continuenewcommand{%
   \@numparams = \ifx\@optionalarg\empty 0\else\@optionalarg \fi \relax
   \@newcommand
}%
%
%
\def\@newcommand#1{%
   \def\@startdef{\expandafter\edef\@commandname}%
   \ifnum\@numparams=0
      \let\@paramdef = \empty
   \else
      \ifnum\@numparams>9
         \errmessage{\the\@numparams\space is too many parameters}%
      \else
         \ifnum\@numparams<0
            \errmessage{\the\@numparams\space is too few parameters}%
         \else
            \edef\@paramdef{%
               \ifcase\@numparams
                  \empty  No arguments.
               \or ####1%
               \or ####1####2%
               \or ####1####2####3%
               \or ####1####2####3####4%
               \or ####1####2####3####4####5%
               \or ####1####2####3####4####5####6%
               \or ####1####2####3####4####5####6####7%
               \or ####1####2####3####4####5####6####7####8%
               \or ####1####2####3####4####5####6####7####8####9%
               \fi
            }%
         \fi
      \fi
   \fi
   \expandafter\@startdef\@paramdef{#1}%
}%
%
%
%
%
\def\@readauxfile{%
   \if@auxfiledone \else 
      \global\@auxfiledonetrue
      \@testfileexistence{aux}%
      \if@fileexists
         \begingroup
            \endlinechar = -1
            \catcode`@ = 11
            \input \jobname.aux
         \endgroup
      \else
         \message{\@undefinedmessage}%
         \global\@citewarningfalse
      \fi
      \immediate\openout\@auxfile = \jobname.aux
   \fi
}%
%
%
\newif\if@auxfiledone
\ifx\noauxfile\@undefined \else \@auxfiledonetrue\fi
%
%
%
%
\@innernewwrite\@auxfile
\def\@writeaux#1{\ifx\noauxfile\@undefined \write\@auxfile{#1}\fi}%
%
%
%
\ifx\@undefinedmessage\@undefined
   \def\@undefinedmessage{No .aux file; I won't give you warnings about
                          undefined citations.}%
\fi
%
%
\@innernewif\if@citewarning
\ifx\noauxfile\@undefined \@citewarningtrue\fi
%
%
%
\catcode`@ = \@oldatcatcode

\def\pfeilso{\leavevmode
            \vrule width 1pt height9pt depth 0pt\relax
           \vrule width 1pt height8.7pt depth 0pt\relax
           \vrule width 1pt height8.3pt depth 0pt\relax
           \vrule width 1pt height8.0pt depth 0pt\relax
           \vrule width 1pt height7.7pt depth 0pt\relax
            \vrule width 1pt height7.3pt depth 0pt\relax
            \vrule width 1pt height7.0pt depth 0pt\relax
            \vrule width 1pt height6.7pt depth 0pt\relax
            \vrule width 1pt height6.3pt depth 0pt\relax
            \vrule width 1pt height6.0pt depth 0pt\relax
            \vrule width 1pt height5.7pt depth 0pt\relax
            \vrule width 1pt height5.3pt depth 0pt\relax
            \vrule width 1pt height5.0pt depth 0pt\relax
            \vrule width 1pt height4.7pt depth 0pt\relax
            \vrule width 1pt height4.3pt depth 0pt\relax
            \vrule width 1pt height4.0pt depth 0pt\relax
            \vrule width 1pt height3.7pt depth 0pt\relax
            \vrule width 1pt height3.3pt depth 0pt\relax
            \vrule width 1pt height3.0pt depth 0pt\relax
            \vrule width 1pt height2.7pt depth 0pt\relax
            \vrule width 1pt height2.3pt depth 0pt\relax
            \vrule width 1pt height2.0pt depth 0pt\relax
            \vrule width 1pt height1.7pt depth 0pt\relax
            \vrule width 1pt height1.3pt depth 0pt\relax
            \vrule width 1pt height1.0pt depth 0pt\relax
            \vrule width 1pt height0.7pt depth 0pt\relax
            \vrule width 1pt height0.3pt depth 0pt\relax}

\def\pfeilsw{ \leavevmode 
            \vrule width 1pt height0.3pt depth 0pt\relax
            \vrule width 1pt height0.7pt depth 0pt\relax
            \vrule width 1pt height1.0pt depth 0pt\relax
            \vrule width 1pt height1.3pt depth 0pt\relax
            \vrule width 1pt height1.7pt depth 0pt\relax
            \vrule width 1pt height2.0pt depth 0pt\relax
            \vrule width 1pt height2.3pt depth 0pt\relax
            \vrule width 1pt height2.7pt depth 0pt\relax
            \vrule width 1pt height3.0pt depth 0pt\relax
            \vrule width 1pt height3.3pt depth 0pt\relax
            \vrule width 1pt height3.7pt depth 0pt\relax
            \vrule width 1pt height4.0pt depth 0pt\relax
            \vrule width 1pt height4.3pt depth 0pt\relax
            \vrule width 1pt height4.7pt depth 0pt\relax
            \vrule width 1pt height5.0pt depth 0pt\relax
            \vrule width 1pt height5.3pt depth 0pt\relax
            \vrule width 1pt height5.7pt depth 0pt\relax
            \vrule width 1pt height6.0pt depth 0pt\relax
            \vrule width 1pt height6.3pt depth 0pt\relax
            \vrule width 1pt height6.7pt depth 0pt\relax
            \vrule width 1pt height7.0pt depth 0pt\relax
            \vrule width 1pt height7.3pt depth 0pt\relax
            \vrule width 1pt height7.7pt depth 0pt\relax
            \vrule width 1pt height8.0pt depth 0pt\relax
            \vrule width 1pt height8.3pt depth 0pt\relax
            \vrule width 1pt height8.7pt depth 0pt\relax
            \vrule width 1pt height9pt depth 0pt\relax
      }


\def\widestnumber#1#2{}

\def\citewarning#1{\ifx\shlhetal\relax 
    \else
    \par{#1}\par
    \fi
}

\def\rm{\fam0 \tenrm}

\def\fakesubhead#1\endsubhead{\bigskip\noindent{\bf#1}\par}



%
%
%

%

\font\textrsfs=rsfs10
\font\scriptrsfs=rsfs7
\font\scriptscriptrsfs=rsfs5

\newfam\rsfsfam
\textfont\rsfsfam=\textrsfs
\scriptfont\rsfsfam=\scriptrsfs
\scriptscriptfont\rsfsfam=\scriptscriptrsfs

\edef\oldcatcodeofat{\the\catcode`\@}
\catcode`\@11

\def\Cal@@#1{\noaccents@ \fam \rsfsfam #1}

\catcode`\@\oldcatcodeofat


\expandafter\ifx \csname margininit\endcsname \relax\else\margininit\fi

\long\def\red#1\endred{}
\long\def\green#1\endgreen{}
\long\def\blue#1\endblue{}
\long\def\private#1\endprivate{}

\def\endred{ \unmatched endred! }
\def\endgreen{ \unmatched endgreen! }
\def\endblue{ \unmatched endblue! }
\def\endprivate{ \unmatched endprivate! }

\ifx\latexcolors\undefinedcs\def\latexcolors{}\fi

\def\emptycs{}
\def\evaluatelatexcolors{%
        \ifx\latexcolors\emptycs\else
        \expandafter\xxevaluate\latexcolors\xxfertig\evaluatelatexcolors\fi}
\def\xxevaluate#1,#2\xxfertig{\setupthiscolor{#1}%
        \def\latexcolors{#2}}


\font\smallfont=cmsl7
\def\rutgerscolor{\ifmmode\else\endgraf\fi\smallfont
\advance\leftskip0.5cm\relax}
\def\setupthiscolor#1{\edef\tmptmpcs{\noexpand\bgroup\noexpand\rutgerscolor
\noexpand\def\noexpand\currentcolor{#1}%
\noexpand}%
\expandafter\let\csname#1\endcsname\tmptmpcs
\def\tmptmpcs{\checkColorUnmatched{#1}\popthecolor}
\expandafter\let\csname end#1\endcsname\tmptmpcs}

\def\checkColorUnmatched#1{\def\expectcolor{#1}%
    \ifx\expectcolor\currentcolor   
    \else \edef\failhere{\noexpand\tryingToClose '\currentcolor' with end\expectcolor}\failhere\fi}

\def\currentcolor{???}

\def\popthecolor{\ifmmode\else\endgraf\fi\egroup}

\expandafter\def\csname#1\endcsname{}

\evaluatelatexcolors

 \let\outerhead\head
 \def\head{\innerhead}
 \let\innerhead\outerhead

 \let\outersubhead\subhead
 \def\subhead{\innersubhead}
 \let\innersubhead\outersubhead

 \let\outersubsubhead\subsubhead
 \def\subsubhead{\innersubsubhead}
 \let\innersubsubhead\outersubsubhead

 \let\outerproclaim\proclaim
 \def\proclaim{\innerproclaim}
 \let\innerproclaim\outerproclaim

 %
 %
 %
 %

\def\demo#1{\medskip\noindent{\it #1.\/}}
\def\enddemo{\smallskip}

\newpage

\head {\S0 Introduction} \endhead  \resetall \sectno=0
\bigskip

Juhasz has asked on the spectrums $c-sp(X) = \{|Y|:Y$ an infinite closed
subspace of $X\}$ and $w-sp(X) = \{w(Y):Y$ a closed subspace of $X\}$.  He
proved \cite{Ju93} that if $X$ is a compact Hausdorff space, then
$|X| > \kappa \Rightarrow c-sp(X) \cap [\kappa,\dsize \sum_{\lambda < \kappa}
2^{2^\lambda}] \ne \emptyset$ and $w(X) > \kappa \Rightarrow w-sp(X) \cap
[\kappa,2^{< \kappa}] \ne \emptyset$.  So under GCH the cardinality spectrum
of a compact Hausdorff space does not omit two successive regular
cardinals, and omit no inaccessible.  Of course, the space $\beta(\omega) 
\backslash
\omega$, the space of nonprincipal ultrafilters on $\omega$, satisfies
$c-sp(X) = \{\beth_2\}$.  Now Juhasz Shelah \cite{JuSh:612} shows that we
can omit many singular cardinals, e.g. under GCH for every regular $\lambda >
\kappa$, there is a compact Hausdorff space $X$ with $c-sp(X) = \{\mu:\mu \le
\lambda,\text{cf}(\mu) \ge \kappa\}$; see more there and in \cite{Sh:652}.
In fact \cite{JuSh:612} constructs a Boolean Algebra, so relevant to the
parallel problems of Monk \cite{M}.  Here we deal with the noncompact case
and get a strong existence theorem.  Note that trivially for a Hausdorff
space $X,|X| \ge \kappa \Rightarrow c-sp(X) \cap [\kappa,2^{2^\kappa}] \ne
\emptyset$, using the closure of any set with $\kappa$ points, so our result
is in this respect best possible.
\bn
We prove
\proclaim{\stag{0.1} Theorem}  For every infinite cardinal $\lambda$ there is
a $T_3$ topological space $X$, even with clopen basis, with $2^{2^\lambda}$
points such that every closed subset with $\ge \lambda$ points has
$|X|$ points.
\endproclaim
\bn
In \S1 we prove a somewhat weaker theorem but with the main points of the
proof present, in \S2 we complete the proof of the full theorem.
\newpage

\head {\S1} \endhead  \resetall \sectno=1
\bigskip

\proclaim{\stag{1.1} Theorem}  Assume $\lambda = \text{\rm cf\/}(\lambda) >
\aleph_0$.  Let $\mu = 2^\lambda,\kappa = \text{\rm Min\/}\{\kappa:2^\kappa >
\mu\}$.  There is a Hausdorff space $X$ with a clopen basis with 
$|X| = 2^\kappa$ such that: if for $Y \subseteq \lambda$ is closed and
$|Y| < |X|$ then $|Y| < \lambda$.
\endproclaim
\bigskip

\demo{Proof}  Let $S \subseteq \{\delta < \kappa:\delta \text{ limit}\}$ be
stationary.  Let $T_\alpha = {}^\alpha \mu$ for $\alpha \le \kappa$ and
let $T = \dbcu_{\alpha \le \kappa} T_\alpha$.  Let $\zeta_\alpha = 
\cup \{\mu \delta + \mu:\delta \in S \cap (\alpha +1)\}$ and let 
$\zeta_{< \alpha} = \cup\{\zeta_\beta:\beta < \alpha\}$.
\enddemo
\bn
\ub{Stage A}:  We shall choose sets $u_\zeta \subseteq T_\kappa$ (for $\zeta <
\mu \times \kappa)$.  Those will be clopen sets generating the topology.  
For each $\zeta$ we choose $(I_\zeta,J_\zeta)$ such that: $I_\zeta$ is a 
$\triangleleft$-antichain of $({}^{\kappa >}\mu,\triangleleft)$ such that 
for every $\rho \in T_\kappa,(\exists ! \alpha)
(\rho \restriction \alpha \in I_\zeta)$ and $J_\zeta \subseteq
I_\zeta$ and we shall let $u_\zeta = \dbcu_{\nu \in J_\zeta}(T_\kappa)
^{[\nu]}$  where $(T_\kappa)^{[\nu]} = \{\rho \in T_\kappa:\nu
\triangleleft \rho\}$.  Let $I_{\alpha,\zeta} = T_\alpha \cap 
I_\zeta,J_{\alpha,\zeta} = T_\alpha \cap J_\zeta$ but we shall have
$\alpha \notin S \Rightarrow I_{\alpha,\zeta} = \emptyset = J_{\alpha,\zeta}$.
\bn
\ub{Stage B}:  Let $Cd:\mu \rightarrow {}^{\lambda^+ >}(T_{< \kappa})$ be onto
such that for every $x \in \text{ Rang}(Cd)$ we have otp$\{\alpha
< \mu:\text{Cd}(\alpha) = x\} = \mu$. \nl
We say $\alpha$ codes $x$ (by Cd) if Cd$(\alpha) = x$.
\bn
\ub{Stage C:Definition}:  For $\delta \le \kappa$ we call $\bar \eta$ a 
$\delta$-candidate if
\mr
\item "{$(a)$}"  $\bar \eta = \langle \eta_i:i \le \lambda \rangle$
\sn
\item "{$(b)$}"  $\eta_i \in T_\delta$
\sn
\item "{$(c)$}"  $(\exists \gamma < \delta)(\dsize \bigwedge_{i < j < \lambda}
\eta_i \restriction \gamma \ne \eta_j \restriction \gamma)$ 
\sn
\item "{$(d)$}"  for every odd $\beta < \delta$, we have \nl
$Cd(\eta_\lambda(\beta)) = \langle \eta_i \restriction \beta:i \le \lambda
\rangle$
\sn
\item "{$(e)$}"  $\eta_\lambda(0)$ codes $\langle \eta_i \restriction
\gamma:i < \lambda \rangle$, where $\gamma = \gamma(\eta \restriction
\lambda) = \text{ Min}\{\gamma < \delta:i < j < \lambda \Rightarrow \eta_i 
\restriction \gamma \ne \eta_j \restriction \gamma\}$, it is well defined 
by clause (c) and
\sn
\item "{$(f)$}"  $\eta_\lambda(0) > \sup\{\eta_i(0):i < \lambda\}$.
\endroster
\bn
\ub{Stage D:Choice}:  Choose $A_{\xi,\varepsilon} \subseteq \lambda$ for 
$\xi < \mu \times \kappa,\varepsilon < \lambda$ such that:

$$
\xi < \mu \times \kappa \and \varepsilon_1 < \varepsilon_2 < \lambda
\Rightarrow |A_{\xi,\varepsilon_1} \cap A_{\xi,\varepsilon_2}| < \lambda
\text{ and even } = \emptyset
$$
\mn
and
$$
\xi_1 < \ldots < \xi_n < \mu \times \kappa,\varepsilon_1 \ldots
\varepsilon_{n_1} < \lambda \Rightarrow \dbca^n_{\ell =1} A_{\xi_\ell,
\varepsilon_\ell} \text{ is a stationary subset of } \lambda.
$$
\mn
Let $\Xi = \bigl\{ \{(\xi_1,\varepsilon_1),\dotsc,(\xi_n,\varepsilon_n)\}:
\xi_1,\dotsc,\xi_n < \mu \times \kappa$ is with no repetitions and
$\varepsilon_1,\dotsc,\varepsilon_n < \lambda \bigr\}$ and for $x \in \Xi$ 
let $A_x = \dbca^n_{\ell=1} A_{\xi_\ell,\varepsilon_\ell}$.
Let $D_0$ be a maximal filter on $\lambda$ 
extending the club filter such that 
$x \in \Xi \Rightarrow A_x \ne \emptyset$ mod $D_0$.
\mn
For $A \subseteq \lambda$ let

$$
{\Cal B}^+(A) = \{x \in \Xi:A \cap A_x = \emptyset \text{ mod } D_0
\text{ but } y \subsetneqq x \Rightarrow A \cap A_y \ne \emptyset
\text{ mod } D_0\}
$$

$$
{\Cal B}(A) =: {\Cal B}^+(A) \cup {\Cal B}^+(\lambda \backslash A).
$$
\bn
\ub{Fact}:  ${\Cal B}(A) =: {\Cal B}^+(A) \cup {\Cal B}^+
(\lambda \backslash A)$ is predense in $\Xi$ i.e.

$$
(\forall x \subseteq \Xi)(\exists y \in {\Cal B}(A))(x \cup y \in \Xi).
$$
\bigskip

\demo{Proof}  If $x \in \Xi$ contradict it then we can add to $D_0$ the set
$\lambda \backslash (A_x \cap A)$ getting $D'_0$.  Now $D'_0$ thus
properly extends $D_0$ otherwise $A_x \cap A = \emptyset$ mod $D_0$
hence, let $x' \subseteq x$ be minimal with this property so $x' \in 
{\Cal B}^+(A)$ and $x$ by assumption satisfies: $\neg(\exists y \in \Xi)
(x \cup y \in {\Cal B}(A))$ so try $y = x$.  For every $z \in \Xi$ we have
$A_z \ne \emptyset$ mod $D_0$.
\enddemo
\bn
\ub{Fact}:  $|{\Cal B}(A)| \le \lambda$ for $A \subseteq \lambda$.
\bigskip

\demo{Proof}   Let ${\bold B}_0$ be the Boolean Algebra
freely generated by $\{x_{\xi,\varepsilon}:\xi < \mu \times \kappa,
\varepsilon < \lambda\}$, by $\Delta$-system argument, except $x_{\xi,
\varepsilon_1} \cap x_{\xi,\varepsilon_2} = 0$ if $\varepsilon_1 \ne
\varepsilon_2$; clearly $\bold B_0$ satisfies $\lambda^+$-c.c.
\sn
Let $\bold B^*$ be the completion of $\bold B_0$.  Let $f^*$ be a 
homomorphism from ${\Cal P}(\lambda)$ into $\bold B^*$ such that 
$C \in D_0 \Rightarrow f^*(C) = 1_{\bold B^*}$ and

$$
f(A_{\xi,\varepsilon}) = x_{\xi,\varepsilon}.
$$
\mn
[Why exists?  Look at the Boolean Algebra ${\Cal P}(\lambda)$ let 
$I_\lambda = \{A \subseteq \lambda:\lambda \backslash A \in D_0\}$ 
and ${\frak A}_0 = I_\lambda \cup
\{\lambda \backslash A:A \in I_\lambda\}$ is a subalgebra of ${\Cal P}
(\lambda)$, and let $I_\lambda \cup \{A_{\xi,\varepsilon}:\xi \le \mu \times
\kappa,\varepsilon = \lambda\}$ generate a subalgebra ${\frak A}$ of
${\Cal P}(\lambda)$; it extends ${\frak A}_0$.  Let 
$f^*_0:{\frak A}_0 \rightarrow \bold B_0$ 
be the homomorphism with kernel $I_\lambda$.  Let $f^*_1$ be the
homomorphism from ${\frak A}$ into $\bold B_0$ extending $f_0$ such that
$f^*_1(A_{\xi,\varepsilon}) = x_{\xi,\varepsilon}$, clearly exists
and is onto.  Now as $\bold B^*$ is a complete Boolean Algebra, $f^*_1$
can be extended to a homomorphism $f^*_2$ from
${\Cal P}(\lambda)$ into $\bold B^*$.  Clearly Ker$(f^*_2) = \text{ Ker}
(f^*_2) = \text{ Ker}(f^*_0) = I_\lambda$ so $f^*_1$ induces an isomorphism
from ${\Cal P}(\lambda)/D_0$ onto Rang$(f^*_1) \subseteq \bold B^*$, so the
problem translates to $\bold B^*$.  So $\bold B_0$ satisfies the
$\lambda^+$-c.c and is a dense subalgebra of $\bold B^*$ hence of 
range$(f^*_2)$, so this range is a $\lambda^+$-c.c. Boolean Algebra hence
${\Cal P}(\lambda)/D_0$ satisfies the fact.] \nl 
Let $\bold B^*_\gamma$ be the complete Boolean subalgebra of 
$\bold B^*$ generated (as a complete subalgebra) by 
$\{x_{\xi,\varepsilon}:\xi < \gamma,\varepsilon < \lambda\}$.  Clearly
$\bold B^* = \dbcu_{\gamma < \kappa} \bold B^*_\gamma$ and $\bold B^*_\gamma$
is increasing with $\gamma$.
\bn
\ub{Stage E}:   We choose by induction on $\delta \in S$ the following
\mr
\item "{$(A)$}"  $w_{\delta,\zeta} \subseteq T_\delta$ (for $\zeta <
\mu \delta + \mu)$ and $J_{\delta,\zeta} \subseteq I_{\delta,\zeta}
\subseteq w_{\delta,\zeta}$
\sn
\item "{$(B)$}"  for each $\delta$-candidate $\bar \eta = \langle
\eta_i:i \le \lambda \rangle$, a uniform filter $D_{\bar \eta}$ on $\lambda$
extending the filter $D_0$
\sn
\item "{$(C)$}"  for each $\nu_1 \ne \nu_2$ in $T_\delta$ for some 
$\zeta < \mu \times \delta + \mu$ we have $\{\nu_1,\nu_2\} \subseteq 
w_{\delta,\zeta}$ and: $(\exists \delta' \in S \cap (\delta + 1))
(\nu_1 \in J_{\delta',\zeta}) \equiv (\exists \delta' \in S \cap
(\delta + 1))(\nu_2 \in J_{\delta',\zeta})$
\sn
\item "{$(D)$}"  if $n < \omega,\mu \times \delta + \mu \le \xi_1 < \ldots 
< \xi_n < \mu \times \kappa$ and $\varepsilon_1,\dotsc,\varepsilon_n < 
\lambda$ then $\dbca^n_{\ell =1} A_{\xi_\ell,\varepsilon_\ell} \ne \emptyset$ 
mod $D_{\bar \eta}$
\sn
\item "{$(E)$}"  if $\delta_1 \in S \cap \delta,\bar \eta$ is a
$\delta$-candidate and $\bar \eta \upharpoonleft \delta_1 = 
\langle \eta_i \restriction \delta_1:i \le \lambda \rangle$ is a
$\delta_1$-candidate \ub{then} 
$D_{\bar \eta \restriction \delta_1} \subseteq D_{\bar \eta}$
\sn
\item "{$(F)_1$}"  $\eta \in w_{\delta,\zeta}$ \ub{iff} $(\exists \delta')
(\delta' \in S \cap (\delta +1) \and \eta \restriction \delta \in
I_{\delta',\zeta})$
\sn
\item "{$(F)_2$}"  if $\bar \eta = \langle \eta_i:i \le \lambda \rangle$ is a
$\delta$-candidate and $\eta_\lambda \in w_{\delta,\zeta}$ \ub{then}
$\{i < \lambda:\eta_i \in w_{\delta,\zeta}\} \in D_{\bar \eta}$ and \nl
$\langle (\exists \delta' \in S \cap (\delta+1))(\eta_\lambda
\restriction \delta' \in J_{\delta',\zeta}) \rangle =$ \nl 
$\text{ LIM}_{D_{\bar \eta}} \langle (\exists \delta' \in S \cap (\delta +1))
(\eta_i \restriction \delta' \in J_{\delta',\zeta}):i < \lambda \rangle$
\sn
\item "{$(F)_3$}"  $w_{\delta,\zeta}$ satisfies the following
{\roster
\itemitem{ $(a)$ }  it is empty if $\zeta < \zeta_{< \delta}$
\sn
\itemitem{ $(b)$ }  has $\le \lambda$ members if $\zeta \in [\zeta_{< \delta},
\zeta_\delta)$
\sn
\itemitem{ $(c)$ }  otherwise $w_{\delta,\zeta}$ is the disjoint union
$w^0_{\delta,\zeta} \cup w^1_{\delta,\zeta} \cup w^2_{\delta,\zeta}$
where \nl
$w^0_{\delta,\zeta} = \bigl \{\eta \in T_\delta:(\exists \delta' \in S 
\cap \delta)(\eta \upharpoonleft \delta' \in w_{\delta',\zeta})\}$ \nl
$w^1_{\delta,\zeta} = \{\eta \in T_\delta:\eta \notin
w^0_{\delta,\zeta}$ and for no $\kappa$-candidate $\bar \eta$ is 
$\eta \triangleleft \eta_\lambda\}$ \nl
$w^2_{\delta,\zeta} = \{\eta \in T_\delta:\eta \notin
w^0_{\delta,\zeta} \cup w^1_{\delta,\zeta}$ and 
for some $\delta$-candidate \nl

$\qquad \qquad \qquad \qquad \bar \eta,\eta_\lambda = \eta$ and
$(\forall i < \lambda)(\exists \delta' \in S \cap \delta)(\eta_i \restriction
\delta' \in w_{\delta',\zeta})$ \nl

$\qquad \qquad \qquad \qquad$ and the set 
$\{i < \lambda:(\exists \delta' \in S \cap \delta)(\eta_i
\restriction \delta' \in J_{\delta,\zeta})\}$ \nl

$\qquad \qquad \qquad \qquad$ or its compliment belongs to 
$D_{{\bar \eta} \restriction \delta^*}$ for some $\delta^* < \delta \bigr\}$
\endroster}
\item "{$(F)_4$}"  $I_{\delta,\zeta} = w^2_{\delta,\zeta} \cup
w^1_{\delta,\zeta}$
\sn
\item "{$(G)$}"  if $\bar \eta$ is a $\delta$-candidate and
$B \subseteq \lambda,f^*(B) \in \bold B^*_{\mu \times (\delta +1)}$,
\ub{then} $B \in D_{\bar \eta} \vee (\lambda \backslash B) \in D_{\bar \eta}$.
\endroster
\noindent
We can ask more explicitly: there is an ultrafilter $D'_{\bar \eta}$ on the
Boolean Algebra $\bold B^*_{\mu \times (\delta +1)}$ such that $D_{\bar \eta} 
= \{B \subseteq \lambda:f^*(B) \in D'_{\bar \eta}\}$. 
\sn
The rest of the proof is split into carrying the construction and proving it
is enough.
\enddemo
\bn
\ub{Stage F:This is Enough}:  First for every $\kappa$-candidate $\bar \eta$
lets $D_{\bar \eta} = \cup\{D_{\bar \nu,\delta}:\delta \in S,\bar \nu$ is a
$\delta$-candidate and $i \le \lambda \Rightarrow \nu_i \triangleleft
\eta_i\}$.  Easily $D_{\bar \eta}$ is a uniform ultrafilter on $\lambda$.  
Let us define the space.
The set of points of the space is $T_\kappa = {}^\kappa \mu$ and a subbase of 
clopen sets will be $u_\zeta$: for $\zeta < \mu \times \kappa$ where
$u_\zeta$ is defined as $u_\zeta =: \cup\{(T_\kappa)^{[\nu]}:\nu 
\in J_\zeta\}$ and 
$J_\zeta =: \dbcu_{\delta \in S} J_{\delta,\zeta}$.  Now note that
\mr
\item "{$(\alpha)$}"  $I_\zeta = \cup\{I_{\delta,\zeta}:\delta \in S\}$ 
is an antichain and $\forall \rho \in T_\kappa \exists
! \delta (\rho \restriction \delta \in I_{\delta,\zeta})$ \nl
[Why?  We prove this by induction on $\rho(0)$ and is straight.  In
details, it is an antichain by the choice $I_{\delta,\zeta} =
w^2_{\delta,\zeta},w^2_{\delta,\zeta} \subseteq T_\delta \backslash
w^0_{\delta,\zeta}$.  As for the second phrase by the first there is at most
one such $\delta$; let $\rho \in T_\kappa$ and assume we have proved it for
every $\rho' \in T_\kappa$ such that $\rho'(0) < \rho(0)$.  By the definition
of $\kappa$-candidate, if there is no $\kappa$-candidate $\bar \eta$ with
$\eta_\lambda = \rho$, \ub{then} for every large enough $\delta \in S$, there
is no $\delta$-candidate $\bar \eta$ with $\eta_\lambda = \rho \restriction
\delta$, hence for any such $\delta,\rho \restriction \delta$ belongs to
$w^0_{\delta,\zeta}$ or to $w^1_{\delta,\zeta}$, in the first case for some
$\delta' \in \delta \cap S$ we have $(\rho \restriction \delta) \restriction
\delta' \in I_{\delta',\zeta}$ so $\rho \restriction \delta' \in
I_{\delta',\zeta}$ and we are done, in the second case $\rho \restriction
\delta \in w^1_{\delta,\zeta} \subseteq I_{\delta,\zeta}$ and we are done.
So assume that there is a $\kappa$-candidate $\bar \eta$ with
$\eta_\lambda = \rho$, by the definition of a candidate it is unique and
$i < \lambda \Rightarrow \eta_i(0) < \rho(0)$, so for each $i < \lambda$
there is $\delta_i \in S$ such that $\eta_i \restriction \delta_i \in
I_{\delta_i,\zeta}$ and let $\gamma = \text{ Min}\{\gamma < \mu:\langle
\eta_i \restriction \gamma:i < \lambda \rangle$ is with no repetition$\}$.  
Let
$A = \{i < \lambda:\eta_i \restriction \delta_i \in J_{\delta,\zeta}\}$ so for
some $\beta < \mu$ we have $f^*_2(A) \in \bold B^*_\beta$.  For $\delta \in
S$, which is $> \sup[\{\gamma,\delta_i:i < \lambda\}]$ we get $\rho
\restriction \delta \in w_{\delta,\zeta}$ and we can finish as before.]
\sn
\item "{$(\beta)$}"  $X$ is a \ub{$T_3$} space \nl
[why?  as we use a clopen basis we really need just 
to separate points which holds by clause (C), i.e. if $\nu_1 \ne \nu_2 \in X$
then for some $\delta \in S$ we have $\nu_1 \restriction \delta \ne \nu_2
\restriction \delta$ and apply clause (C) to $\nu_1 \restriction \delta,
\nu_2 \restriction \delta$]
\sn
\item "{$(\gamma)$}"  $|X| = \mu^\kappa = 2^\kappa$ \nl
[why?   as $T_\kappa$ is the set of points of $X$]
\sn
\item "{$(\delta)$}"  suppose $Y = \{\eta_i:i < \lambda\} \subseteq X =
T_\kappa$ and $\dsize \bigwedge_{i < j} \eta_i \ne \eta_j$.  We need to
show that $|c \ell(Y)|$ large, i.e. has cardinality $2^\kappa$.
\ermn
Choose $\gamma$ such that $\langle \eta_i \restriction \gamma:i < \lambda
\rangle$ is with no repetitions.
\mn
Let 

$$
\align
W_{\bar \eta} = \{<>\} \cup \bigl\{\rho:&\text{for some } \alpha \le \kappa,
\rho \in T_\alpha,\rho(0) \text{ code } \langle \eta_i \restriction
\gamma:i < \lambda \rangle, \\
   &\rho(0) > \sup\{\eta_i(0):i < \lambda\} \text{ and} \\
  &(\forall \beta < \ell g(\rho))(\beta \text{ odd } \Rightarrow \rho(\beta)
\text{ code } \langle \eta_i \restriction \beta:i < \lambda \rangle \char 94
\langle \rho \restriction \beta \rangle) \bigr\}.
\endalign
$$
\mn
So clearly:
\mr
\widestnumber\item{$(iii)$}
\item "{$(i)$}"   $W_{\bar \eta} \cap T_1 \ne \emptyset$
\sn
\item "{$(ii)$}"  $W_{\bar \eta}$ is a subtree of 
$(\dbcu_{\alpha \le \kappa} T_\alpha,\triangleleft)$
(i.e. closed under initial segments, closed under limits),
\sn
\item "{$(iii)$}"   every $\rho \in W_{\bar \eta} \cap T_\alpha$ where
$\alpha < \kappa$ has a successor and if $\alpha$ is even has $\mu$ 
successors.
\ermn
So $|W_{\bar \eta} \cap T_\kappa| = \mu^\kappa$.
\mn
So enough to prove
\mr
\item "{$(*)$}"  if $\rho \in W_{\bar \eta} \cap T_\kappa$ then
$\rho \in c \ell\{\eta_i:i < \lambda\}$.
\ermn
Let $\bar \eta = \langle \eta_i:i < \lambda \rangle,\eta_\lambda = \rho,
\bar \eta' = \bar \eta \char 94 \langle \rho \rangle$ and the filter
$D_{{\bar \eta}'} = \cup\{D_{\langle \bar \eta'_i \upharpoonleft \delta:
i \le \lambda \rangle}:\delta \in S \text{ and } \delta \ge \gamma\}$
is a filter by clause (E) and even ultrafilter by clause (G).
\mn
Now for every $\zeta$, by clause (F)$_2$ for $\delta$ large enough

$$
\text{Truth Value}(\rho \in u_\zeta) = 
\text{ lim}_{D_{\langle \bar \eta'_i \upharpoonleft \delta:
i \le \delta \rangle}}
\langle \text{Truth Value}(\eta_i \in u_\zeta):i < \lambda \rangle.
$$
As $\{u_\zeta:\zeta < \mu \times \kappa\}$ 
is a clopen basis of the topology, we are done.
\bn
\ub{Stage G: The construction}:

We arrive to stage $\delta \in S$.  So for every $\delta$-candidate
$\bar \eta = \langle \eta_i:i \le \lambda \rangle$, let

$$
D'_{\bar \eta} = \cup\{D_{\langle \eta_i \upharpoonleft \delta_1:
i \le \lambda \rangle}:\delta_1 \in \delta \cap S \text{ and } 
\langle \eta_i \upharpoonleft \delta_1:i \le \lambda \rangle 
\text{ a } \delta_1 \text{-candidate}\} \cup D_0.
$$
\mn
\ub{Note}:  $|T_\delta| = \mu$ by the choice of $\kappa$.

Let $<^*_\delta$ be a well ordering of $T_\delta$ such that:
$\nu_1(0) < \nu_2(0) \Rightarrow \nu_1 <^*_\delta \nu_2$.
\sn
Hence
\mr
\item "{$(*)$}"  $\langle \eta_i:i \le \lambda \rangle$ a $\delta$-candidate
$\Rightarrow \dsize \bigwedge_{i < \lambda} \eta_i <^*_\delta \eta_\lambda$.
\ermn
So let $\{\langle \nu_{1,\zeta},\nu_{2,\zeta} \rangle:\zeta_{< \delta} \le
\zeta < \zeta_\delta\}$ list $\{(\nu_1,\nu_2):\nu_1 <^*_\delta \nu_2\}$;
such a list
exists as $\zeta_\delta \ge \zeta_{< \delta} + \mu$ and $|T_\delta| = \mu$.  
Now we choose by induction on $\zeta < \zeta_\delta$ the following
\mr
\item "{$(\alpha)$}"  $D^\zeta_{\bar \eta}$ for $\bar \eta$ a
$\delta$-candidate when $\zeta \ge \zeta_{< \delta}$
\sn
\item "{$(\beta)$}"  $w^*_{\delta,\zeta},I_{\delta,\zeta},J_{\delta,\zeta}$
\sn
\item "{$(\gamma)$}"  $D^{\zeta_{< \delta}}_{\bar \eta}$ is
$D'_{\bar \eta}$ which was defined above
\ermn
such that
\mr
\item "{$(\delta)$}"  $D^\zeta_{\bar \eta}$ for 
$\zeta$ in $[\zeta_{< \delta},\zeta_\delta]$ is increasing 
continuous
\sn
\item "{$(\varepsilon)$}"  if $n < \omega,\zeta_{< \delta} \le 
\zeta \le \xi_1 < \xi_2 < \ldots < \xi_n < \mu \times \kappa$ and 
$\varepsilon_1,\dotsc,\varepsilon_n < \lambda^+$ then $\dbca^n_{\ell =1} 
A_{\xi_\ell,\varepsilon_\ell} \ne \emptyset$ mod $D^\zeta_{\bar \eta}$
\sn
\item "{$(\zeta)$}"  
$D^{\zeta +1}_{\bar \eta},I_{\delta,\zeta},J_{\delta,\zeta}$ satisfies the
requirement (F)$_2$
\sn
\item "{$(\eta)$}"  $\nu_{1,\zeta} \in J_{\delta,\zeta} \Leftrightarrow
\nu_{2,\zeta} \notin J_{\delta,\zeta}$ \ub{or} $\nu_{1,\zeta},\nu_{2,\zeta}
\in w^0_{\delta,\zeta}$
\sn
\item "{$(\theta)$}"  $D^\zeta_{\bar \eta}$ is $D'_{\bar \eta}
+ \{A_{\zeta_1,\varepsilon_{\bar \eta}(\zeta_0)}:\zeta_1 < \zeta\}$ for some
function 
$\varepsilon_{\bar \eta}:[\zeta_{< \delta},\zeta) \Rightarrow \lambda$.
\ermn
\ub{Note}:  For $\zeta=0$, condition $(\varepsilon)$ holds by the induction
hypothesis (i.e. clause (D)) and choice of $D'_\eta$ (and choice of the
$A_{\xi,\varepsilon}$'s if for no $\delta_1,\bar \eta \restriction \delta_1$
is a $\delta_1$-candidate).
\mr
\item "{$(\iota)$}"  if $\zeta < \zeta_{< \delta}$ \ub{then}:
\ermn

$$
w_{\delta,\zeta} = w^0_{\delta,\zeta} \cup w^1_{\delta,\zeta} \cup
w^2_{\delta,\zeta} \text{ are defined as in } (F)_2
$$

$$
I^\zeta_{\delta,\zeta} = w^1_{\delta,\zeta} \cup w^2_{\delta,\zeta}
$$

$$
\align
J^\zeta_{\delta,\zeta} = \{\eta \in T_\delta:&\delta \in w^2_{\delta,\zeta}
\text{ and for some } 
\delta \text{-candidate } \bar \eta \text{ we have } \eta_\lambda = \eta \\
  &\text{hence }(\forall i < \lambda)(\exists \delta' \in S \cap \delta)
[\eta_i \restriction \delta' \in w_{\delta',\zeta}] \\
  &\text{and } \{i < \lambda:(\exists \delta' \in S \cap \delta)
[\eta_i \restriction \delta' \in J_{\delta',\zeta}]\} \text{ belongs to } 
D'_{\bar \eta}\}.
\endalign
$$
\mn
[Note in the context above, by the induction hypothesis $(\exists \delta' \in
S \cap \delta)[\eta_i \restriction \delta' \in w_{\delta',\zeta}]$ is
equivalent to $(\exists \delta' \in S \cap \delta)[\eta_i \restriction
\delta' \in I_{\delta',\zeta}]$ and thus $\delta'$ is unique.  Of course,
they have to satisfy the relevant requirements from (A)-(G)].
\mn
The cases $\zeta \le \zeta_{< \delta},\zeta$ limit are easy. \nl
The crucial point is: we have $\langle D^\zeta_{\bar \eta}:\bar \eta$ a
$\delta$-candidate$\rangle$ and $\zeta \in [\zeta_{< \delta},\zeta_\delta)$ 
and we should define $w_{\delta,\zeta},I_{\delta,\zeta}$ and
$D^{\zeta +1}_{\bar \eta}$ to which the last stage is dedicated. 
\bn
\ub{Stage H}:  Define by induction on $n < \omega$,

$$
w^\zeta_0 = \{\nu_{1,\zeta},\nu_{2,\zeta}\}
$$

$$
w^\zeta_{n+1} = \{\eta^\rho_i:i < \lambda,\rho \in w_n \text{ and }
\bar \eta^\rho \text{ is a } \delta\text{-candidate with } \eta^\rho_\lambda
= \rho\}.
$$
\mn
Note that $\eta^\rho_i <^*_\delta \rho$. \nl
Let $w = w_{\delta,\zeta} = I_{\delta,\zeta} = \dbcu_{n < \omega} w^\zeta_n$,
so $|w_{\delta,\zeta}| \le \lambda$ (note that this is the first ``time"
we deal with $\zeta$).
\mn
We need: to choose $J_{\alpha,\zeta} \cap w_{\delta,\zeta}$ so that 
the cases of condition $(\zeta)$ (i.e. (F)$_2$) for
$\bar \eta^\rho,\rho \in w$ hold and condition $(\eta)$ (i.e. (C) for
$\nu_{1,\zeta},\nu_{2,\zeta}$) holds.
\sn
Let $w'_{\delta,\zeta} = \{\rho \in w_{\delta,\zeta}:\bar \eta^\rho$ is 
well defined$\}$, (so $w'_{\delta,\zeta} \subseteq w_{\delta,\zeta}$).  
Let $w'_{\delta,\zeta} = \{\rho[\zeta,\varepsilon]:\varepsilon < \varepsilon^*
\le \lambda\}$.  Now we define 
$D^{\zeta +1}_{\bar \eta^{\rho[\zeta,\varepsilon]}}$
as $D^\zeta_{\eta^{\rho[\zeta,\varepsilon]}} + A_{\zeta,\varepsilon}$,
clearly ``legal".
\sn
Let $A'_{\zeta,\varepsilon} = \{i < \lambda:i \in A_{\zeta,\varepsilon}$ and
$i > \varepsilon$ and $\eta^{\rho[\zeta,\varepsilon]}_i \notin
\{\eta^{\rho[\zeta,\varepsilon_1]}_{i_1}:\varepsilon_1 < i$ and $i_1 < i\}$ 
and $\eta^{\rho[\zeta,\varepsilon]}_i \ne \nu_{1,\zeta},\nu_{2,\zeta}\}$.
\nl
Observe
\mr
\item "{$(*)_1$}"   $A_{\zeta,\varepsilon} \backslash A'_\varepsilon$ is not
stationary by Fodor's lemma as $\langle \eta^{\rho[\varepsilon]}_i:i <
\lambda \rangle$ is with no repetition.
\ermn
Now we shall prove that
\mr
\item "{$(*)_2$}"  the sets $\{\eta^{\rho[\varepsilon]}_i:
i \in A'_\varepsilon\}$ for $\varepsilon > \varepsilon^*$ 
are pairwise disjoint.
\ermn
So toward contradiction suppose 
$i_1 \in A'_{\varepsilon_1},i_2 \in A'_{\varepsilon_2},
\varepsilon_1 < \varepsilon_2 < \varepsilon^*$ and 
$\eta^{\rho^{[\zeta,\varepsilon_1]}}_{i_1} = 
\eta^{\rho^{[\zeta,\varepsilon_2]}}_{i_2}$
and try to get a contradiction.
\bn
\ub{Case 1}:  $i_2 > i_1$.  

As $i_1 \in A'_{\varepsilon_1}$ we have $i_1 > \varepsilon_1$ 
similarly $i_2 > \varepsilon_2$ but $\varepsilon_1 <
\varepsilon_2$ so $i_2 > \varepsilon_2 > \varepsilon_1$, and by the
assumption $i_2 > i_1$.  So $\eta^{\rho^{[\zeta,\varepsilon_1]}}_{i_1}$ 
belongs to the set $\{\eta^{\rho^{[\zeta,\varepsilon]}}_i:\varepsilon < i_2 
\and i < i_2\}$ so $\eta^{\rho[\zeta,\varepsilon_2]}_{i_2} \ne 
\eta^{\rho[\zeta,\varepsilon_1]}_{i_1}$
as $\eta^{\rho[\zeta,\varepsilon_2]}_{i_2}$ does not belong to this set as 
$i_2 \in A'_{\varepsilon_2}$.
\bn
\ub{Case 2}:  $i_2 < i_1$.

As $i_2 \in A'_{\zeta,\varepsilon_2}$ necessarily $\varepsilon_2 < i_2$.
So $\varepsilon_2 < i_2 < i_1$ so $\eta^{\rho^{[\zeta,\varepsilon_2]}}_{i_2} 
\in \{\eta^{\rho^{[\varepsilon]}}_i:\varepsilon < i_1 \and \ell^i < i_1\}$
but $\eta^{\rho^{[\zeta,\varepsilon_1]}}_{i_2}$ does not belong to this set as
$i_1 \in A'_{\varepsilon_1}$ hence $\eta^{[\zeta,\varepsilon_1]}_{i_1},
\eta^{[\zeta,\varepsilon_2]}_{i_2}$ cannot be equal.
\bn
\ub{Case 3}:  $i_1 = i_2$.  

As $i_1 \in A'_{\varepsilon_1}$ we have
$i_1 \in A_{\zeta,\varepsilon_1}$ similarly $i_2 \in A_{\zeta,\varepsilon_2}$
but those sets are disjoint; a contradiction. \nl
So $(*)_2$ holds.
\mn
Now define $w^{\zeta,\ell}_n$ for $\ell = 1,2,n < \omega$ by induction on

$$
n:w^{\zeta,\ell}_0 = \{\nu_{\ell,\zeta}\}
$$

$$
w^{\zeta,\ell}_{n+1} = \{\eta^{\rho^{[\zeta,\varepsilon]}}_i:
\rho[\zeta,\varepsilon] \in w^{\zeta,\ell}_n \text{ and } 
i \in A'_\varepsilon \text{ and } \varepsilon < \varepsilon^*\}.
$$
\mn
Let $w^{\zeta,\ell} = \dbcu_{n < \omega} w^{\zeta,\ell}_n$, now by 
$(*)_2$, $w^{\zeta,1} \cap w^{\zeta,2} = \emptyset$
(note the clause $\eta^{\rho^{[\zeta,\varepsilon]}}_i 
\ne \nu_{1,\zeta}$ in the definition of $A'_\varepsilon$). \nl
So we define

$$
J_{\delta,\zeta} = w^{\zeta,2}.
$$
\mn
Now it is easy to check clause (F), i.e. $(\zeta)$ and 
we have finished the induction on $\zeta < \zeta_\delta$.  Now choose
$D_{\bar \eta}$ to satisfy clause (G) and to extend $\dbcu_{\zeta < \zeta
_\delta} D^\zeta_{\bar \eta}$, so we are done. \nl
${{}}$  \hfill$\square_{\scite{1.1}}$
\bn
\centerline{$* \qquad * \qquad *$}
\newpage

\head {\S2 The singular case and the full result} \endhead  \resetall \sectno=2
\bigskip

\proclaim{\stag{2.1} Theorem}  Assume $\lambda > \text{\rm cf\/}(\lambda)$.  
Let $\mu = 2^\lambda,\kappa = \text{\rm Min\/}\{\kappa:2^\kappa >
\mu\}$.  There is a Hausdorff space $X$ with a clopen basis with 
$|X| = 2^\kappa$ such that for $Y \subseteq \lambda$ closed $|Y| < |X|
\Rightarrow |Y| < \lambda$.
\endproclaim
\bigskip

\demo{Proof}  For $\lambda$ singular we should replace the filter 
$D_0$ on $\lambda$.  So let $\lambda = \dsize \sum_{j <
\text{cf}(\lambda)} \lambda_j,\lambda_j$ strictly increasing $\bar \lambda =
\langle \lambda_j:j < \text{cf}(\lambda) \rangle$.  Let $D^*_{\bar \lambda}
= \{A \subseteq \lambda:\text{for every } j < \text{ cf}(\lambda) \text{ large
enough, the set } 
A \cap \lambda^+_j \text{ contains a club of } \lambda^+_j\}$.
\mn
We can find a partition $\langle A^j_\alpha:\alpha < \lambda^+_j \rangle$
of $\lambda^+_j \backslash \lambda_j$ to stationary sets; let us
stipulate $A^j_\alpha
= \emptyset$ when $\lambda^+_j \le \alpha < \lambda$ and let $\bar A^* =
\langle A_\alpha = \dbcu_{j < \text{ cf}(\lambda)} A^j_\alpha:\alpha < \lambda
\rangle$ (so $A_\alpha \ne \emptyset \text{ mod } D^*_\lambda$ 
and $\alpha < \beta < \lambda
\Rightarrow A_\alpha \cap A_\lambda = \emptyset$).  Let $\{f_\xi:\xi < \mu
\times \kappa\}$ be a family of functions from $\lambda$ to $\lambda$ such
that if $n < \omega,\xi_1 < \ldots < \xi_n < \mu \times \kappa$ and
$\varepsilon_1,
\dotsc,\varepsilon_n < \lambda$ then $\{\alpha < \lambda:f_{\varepsilon_\ell}
(\alpha) = \varepsilon_\ell$ for $\ell = 1,\dotsc,n\}$ is not empty (exists
by \cite{EK}).  Now for $\xi < \mu \times \kappa$ and $\varepsilon < \lambda$
we let $A_{\xi,\varepsilon} = \cup \{A_\alpha:f_\xi(\alpha) = \varepsilon\}$.
Clearly $\xi < \mu \times \kappa \and \varepsilon_1 < \varepsilon_2 <
\lambda \Rightarrow A_{\xi,\varepsilon_1} \cap A_{\xi,\varepsilon_2} =
\emptyset$, and also: if $n < \omega,\xi_1 < \ldots < \xi_n < \mu \times
\kappa$ and $\varepsilon_1,\dotsc,\varepsilon_n < \lambda$ then
$\dbca^n_{\ell =1} A_{\xi_\ell,\varepsilon_\ell} \ne \emptyset$ mod 
$D^*_\lambda$.  Let $D_0$ be a maximal 
filter on $\lambda$ extending $D^*_\lambda$
and still satisfying $\dbca^n_{\ell =1} A_{\xi_\ell,\varepsilon_\ell} \ne
\emptyset$ mod $D_0$ for $n,\xi_\ell,\varepsilon_\ell (\ell < n)$ as above.
\mn
Now the proof proceeds as before.   All is the same except in stage H where
we use $\lambda$ regular, $D_0$ contains all clubs of $\lambda$.
\mn
The point is that we define $A'_\varepsilon$ as before, the main question is:
why $A'_\varepsilon = A_\varepsilon$ mod $D^*_{\bar \lambda}$.
\sn
Choose $j^* < \text{ cf}(\lambda)$ such that:

$$
\varepsilon < \lambda_{j^*}.
$$
\mn
So it is enough to show
\mr
\item "{$(*)$}"  if $j^* \le j < \text{ cf}(\lambda)$ then \nl
$A'_\varepsilon \cap [\lambda_j,\lambda^+_j) = A_\varepsilon \cap [\lambda_j,
\lambda^+_j)$ mod $D_{\lambda^+_j}$
\ermn
(where $D_{\lambda^+_j}$-the club filter on $\lambda^+_j$). \nl
Looking at the definition of $A'_{\zeta,\varepsilon}$,

$$
\align
A'_{\zeta,\varepsilon} \cap [\lambda_j,\lambda^+_j) = \bigl\{i \in [\lambda_j,
\lambda^+_j):&i \in A_{\zeta,\varepsilon} \cap [\lambda_j,\lambda^+_j) \\
  &\text{and }\eta^{\rho[\zeta,\varepsilon]}_{i_1} \notin \{
\eta^{\rho[\zeta,\varepsilon_1]}_{i_1}:\varepsilon_1 < i \text{ and} \\
  &i_1 < i\} \text{ and }
\eta^{\rho[\varepsilon]}_i \ne \nu_{1,\zeta} \bigr\}
\endalign
$$
\mn
as $\langle \eta^{\rho^{[\zeta,\varepsilon]}}_i:\lambda_j \le i < 
\lambda^+_j \rangle$
is with no repetition and Fodor's theorem holds (can formulate the demand
on $D$).  Just check that the use of $A'_{\zeta,\varepsilon}$ in \S1 still
works.
\enddemo
\bn
\margintag{2.2}\ub{\stag{2.2} Conclusion}:  If $\lambda \ge \aleph_0,\kappa = \text{ Min}
\{\kappa:2^\kappa > 2^\lambda\}$, \ub{then} there is a $T_3$-space
$\lambda,|X| = 2^\kappa$ with no closed subspace of cardinality $\in
[\lambda,2^\kappa)$.  \hfill$\square_{\scite{2.1}}$
\bn
\centerline{$* \qquad * \qquad *$}
\bn
We still would like to replace $2^\kappa$ by $2^{2^\lambda}$.
\bigskip

\proclaim{\stag{2.3} Theorem}  For $\lambda \ge \aleph_0$ there is a $T_3$
space $X$ with clopen basis such that: no closed subspace has cardinality in
$[\lambda,2^{2^\lambda}]$.
\endproclaim
\bigskip

\demo{Proof}   For $\lambda = \aleph_0$ it is known so let $\lambda >
\aleph_0$.  Like the proof of \scite{1.1} with $\kappa = 2^\mu$.
\sn
The only problem is that 
$T_\delta = {}^\delta \mu$ may have cardinality $> 2^\mu$ so we have to
redefine a $\delta$-candidate (as there are too many $\eta_i \restriction
\gamma$ to code) and in the crucial Stages G and H we have the list 
$\{(\nu^\delta_{1,\varepsilon},
\nu^\delta_{2,\varepsilon}):\varepsilon < |T_\delta|\}$ but possibly
$|T_\delta| > 2^\mu$.  Still $|T_\delta| \le \mu^{|\delta|} \le 2^\mu$; so
instead dedicating one $\zeta \in [\zeta_{< \delta},\zeta_\delta)$ to deal
with any such pair we just do it for each ``kind" of pairs such that the
number of kinds is $\le \mu$, (but we can deal with all of them at once).
\bn
\ub{Stage $B'$}:

Let $Cd:\mu \rightarrow {\Cal H}_{< \lambda^+}(\mu)$ be such that for every
$x \in {\Cal H}_{< \lambda^+}(\mu)$ for $\mu$ ordinals $\alpha < \mu$ we have
$Cd(\alpha) = x$.
\bn
\ub{Stage $C'$}:

For limit $\delta \le \kappa$ we call $\bar \eta$ a $\delta$-candidate if:
\mr
\item "{$(a)$}"  $\bar \eta = \langle \eta_i:i \le \lambda \rangle$
\sn
\item "{$(b)$}"  $\eta_i \in T_\delta$
\sn
\item "{$(c)$}"  for some $\gamma,\langle \eta_i \restriction \gamma:i <
\lambda \rangle$ is with no repetition
\sn
\item "{$(d)$}" for odd $\beta < \delta$ we have \nl
$Cd(\eta_\lambda(\beta)) = \langle(\eta_i(\beta -1),\eta_i(\beta)):i <
\lambda \rangle$
\sn
\item "{$(e)$}"  $Cd(\eta_\lambda(0)) = \{(i,j,\gamma,\eta_i(\gamma),
\eta_j(\gamma)):i < j < \lambda \text{ and for some } i_1 < j_1 < \lambda,
\gamma \text{ minimal such that } \eta_{i_1}(\gamma) \ne \eta_{j_1}(\gamma)\}$
\sn
\item "{$(f)$}"  $\eta_\lambda(0) > \sup\{\eta_i(0):i < \lambda\}$.
\ermn
So
\mr
\item "{$(*)_1$}"  if $\langle \eta_i:i \le \lambda \rangle$ is a
$\delta_1$-candidate, $\delta_0 < \delta_1$ limit and
$(\exists \gamma < \delta_0)(\langle \eta_i \restriction \gamma:
i \le \lambda \rangle$ with no repetitions \ub{then} $\langle \eta_i
\restriction \delta_0:i \le \lambda \rangle$ is a $\delta_0$-candidate
\sn
\item "{$(*)_2$}"  if $\eta_i \in T_\kappa$ for $i < \kappa$ are pairwise
distinct \ub{then} for $2^\mu$ sequences 
$\eta_\lambda \in T_\kappa$ we have $\langle \eta_i:i \le \lambda \rangle$ 
is a $\kappa$-candidate.
\endroster
\bn
\ub{Stage H'}:

For each $\varepsilon < |T_\delta|$ we can choose $v_{\delta,\varepsilon} =
\cup\{v_{\delta,\varepsilon,n}:n < \omega\}$ where we define
$v_{\delta,\varepsilon,n}$ by induction on $n$ as follows:
\sn
$v_{\delta,\varepsilon,0} = \{\nu^\delta_{1,\varepsilon},
\nu^\delta_{2,\varepsilon}\},v_{\delta,\varepsilon,n+1} = v_{\delta,
\varepsilon,n} \cup \{\eta^\rho_i:\rho \in v_{\delta,\varepsilon,n}$ and
$\bar \eta^\rho$ is a $\delta$-candidate such that $\eta^\rho_\lambda =
\rho\}$.
We choose $u_\varepsilon = u_{\delta,\varepsilon} \in 
[\delta]^{\le \lambda}$ such
that: if $\bar \eta$ is a $\delta$-candidate satisfying $\eta_\lambda \in
v_{\delta,\varepsilon}$ (so $\eta_i \in v_{\delta,\varepsilon}$ for
$i < \lambda$) then $0 \in u_\varepsilon \and i < j < \lambda \Rightarrow
\text{ Min}\{\gamma:\eta_i(\gamma) \ne \eta_j\gamma)\} \in u_\varepsilon$.

As $|T_\delta| \le 2^\mu$ and $\mu^\lambda = \mu$ by 
Engelking Karlowic \cite{EK} there are functions $H^\delta_\Upsilon:
T_\delta \rightarrow {\Cal H}_{< \lambda^+}(\mu)$ for 
$\Upsilon \in [\zeta_{< \delta},
\zeta_\delta)$ such that for every $w \in [T_\delta]^\lambda$ and
$h:w \rightarrow {\Cal H}_{< \lambda^+}(\mu)$ there is 
$\Upsilon \in [\zeta_{< \delta},
\zeta_\delta)$ such that $h \subseteq H^\delta$.
\mn
As $\mu = \mu^\lambda = |{\Cal H}_{< \lambda^+}(\mu)|$, \wilog \, 
$|\text{Rang}(H^\delta_\Upsilon)| \le
\lambda$ (divide $H^\delta_\Upsilon$ to $\le 2^\lambda = \mu$ functions).
\mn
For each $\varepsilon < |T_\delta|$ let $h^\varepsilon_\delta:
v_{\delta,\varepsilon} \rightarrow {\Cal H}_{< \lambda^+}(\mu)$ be
$h^\varepsilon_\delta(\eta) = (h^{\varepsilon,0}_\delta(\eta),
h^{\varepsilon,1}_\delta(\eta),h^{\varepsilon,2}_\delta(\eta))$ where

$$
h^{\varepsilon,0}_\delta(\eta) = \text{ otp}(\{\nu \in w^\varepsilon_\delta:
\nu <^*_\delta \eta\},<^*_\delta)
$$

$$
h^{\varepsilon,1}_\delta(\eta) = \{\langle \gamma,\eta(\gamma) \rangle:
\gamma \in u_{\delta,\varepsilon}\}
$$

$$
h^{\varepsilon,2}_\delta(\eta) = \text{ truth value of }
\eta \in v_{\delta,\varepsilon,0}
$$
\mn
(the function $h^\varepsilon_\delta$ belongs to 
${\Cal H}_{< \lambda^+}(\mu)$ as $|v_{\delta,\varepsilon}| \le \lambda$); let

$$
\Upsilon_\varepsilon = \text{ Min}\{\Upsilon \in [\zeta_{< \delta},
\zeta_\delta):h^\varepsilon_\delta \subseteq H^\delta_\Upsilon\}
$$
\mn
(well defined).  Let $\gamma^\delta_\Upsilon =: 
\sup\{\gamma < \lambda^+:\gamma$ is the first cardinal in some sequence
$\bar \lambda \text{ from } \text{(Rang}(H^\delta_\Upsilon)\}$, let 
$g^\delta_\Upsilon$ be a one-to-one function from
$\gamma^\delta_\Upsilon$ into $\lambda$.
\mn
Next we can define the $D^\Upsilon_{\bar \eta}$ for $\bar \eta$ a
$\delta$-candidate; for $\Upsilon < \mu$:

$$
D^{\Upsilon +1}_{\bar \eta} = D^\Upsilon_{\bar \eta} + A_{\Upsilon,
\gamma^\delta_\Upsilon}.
$$
\mn
In Stage $\Upsilon \in [\zeta_{< \delta},\zeta_\delta)$ we deal with all
$\varepsilon < |T_\delta|$ such that $\Upsilon_\varepsilon = \Upsilon$.  Now
we treat the choice of $I_{\delta,\zeta},J_{\delta,\zeta},w_{\delta,\zeta}$.
We can finish as before (but dealing with many cases at once).
\hfill$\square_{\scite{2.3}}$
\enddemo

\newpage
    
REFERENCES.  
\bibliographystyle{lit-plain}
\bibliography{lista,listb,listx,listf,liste}

\shlhetal
\enddocument

\bye